\documentclass[aps,prl,twocolumn,superscriptaddress,runinaddress,floatfix,preprintnumbers,showpacs,nofootinbib]{revtex4-1}

\usepackage{dcolumn}
\usepackage{amssymb}
\usepackage{amsmath}
\usepackage{amsthm}
\usepackage{latexsym}
\usepackage{graphicx}
\usepackage{color}
\usepackage{ulem}
\usepackage{braket}

\newcommand*{\xdash}[1][3em]{\rule[0.5ex]{#1}{0.55pt}}

\begin{document}

\preprint{}

\title{In principle determination of generic priors}
\author{Cael L. Hasse}
\email{Electronic address: cael.hasse@adelaide.edu.au}
\affiliation{Special Research Centre for the Subatomic Structure of Matter and Department of Physics, University of Adelaide 5005, Australia.}
\date{\today}
\begin{abstract}
Probability theory as extended logic is completed such that essentially any probability may be determined. This is done by considering propositional logic (as opposed to predicate logic) as syntactically sufficient and imposing a symmetry from propositional logic. It is shown how the notions of `possibility' and `property' may be sufficiently represented in propositional logic such that 1) the principle of indifference drops out and becomes essentially combinatoric in nature and 2) one may appropriately represent assumptions where one assumes there is a space of possibilities but does not assume the size of the space.
\end{abstract}

\maketitle

\section{Introduction}

This article is a summation of current and ongoing research. Some work in the literature may not yet be properly considered. \\

It can be argued that Bayesian probability theory is the general framework for scientific prediction and inference \cite{Jaynes:2003, Finetti}. It is a calculus for normative statements
\begin{equation}
P(A|B), \nonumber
\end{equation}
where $A$ and $B$ are propositions. These statements encode degrees of belief/plausibility of $A$, given $B$. Fully half of probability theory is encoded into two simple rules\footnote{The notation we are using corresponds to $AB=A$ and $B, A+B=A$ or $B$, and $\bar{A}=$ not $A$},
\begin{align*}
  \textrm{Product rule:} \hspace{2mm}& P(AB|C)=P(A|C)P(B|AC) \\
&\hspace{14.7mm}=P(B|C)P(A|BC) \\
\textrm{Sum rule:} \hspace{2mm}& P(A|B)+P(\bar{A}|B)=1,
\end{align*}
from which we have a generalised sum rule
\begin{equation*}
P(A+B|C) = P(A|C) + P(B|C) - P(AB|C).
\end{equation*}
These rules give relationships between different probabilities but do not constrain the probabilities enough to uniquely determine them  \cite{Cox:1946}. This is the incompleteness of probability theory.\\

In particular, in inference we often want to determine the probability of a hypothesis $H$ given data $D$ and perhaps some background knowledge $I$,
\begin{equation}
P(H|DI)=\frac{P(D|HI)P(H|I)}{P(D|I)}. \nonumber
\end{equation}
Probabilities like $P(H|I)$ are often called `priors'. These need to be determined to work out $P(H|DI)$.\\

Many methods have been invented to determine priors for different situations. If one is a Subjective Bayesian where a probability is a degree of belief relative to some agent, the (perhaps counterfactual) agent is free to just choose the value of their priors based on intuition or gut instinct. If one is an Objective Bayesian, one wants to find methods of derivation such that each probability assignation can be considered unique and agent independent (although different agents may make different assumptions and so may still consider different degrees of plausibility; i.e., Objective Bayesians still consider probability theory subjective in one sense). Such methods include Laplace's principle of indifference, transformation group methods, and maximum entropy methods \cite{Shore:1980}. However, methods such as these are not always of use for calculation of generic priors.\\

I propose a method whereby in principle any prior may be determined. Moreover, it is a method for calculating essentially\footnote{We shall be using a finite sets policy where one always starts with finite sets and only then takes limits. See \cite{Jaynes:2003} for detailed motivation.} any\footnote{There is also no general method for calculating factors $P(D|HI)$ - called likelihoods - unless $H$ predicts $D$ or $\bar{D}$ for certain; the calculation often reduces to calculation of a prior.} probability. This will be accomplished by completing the Objective Bayesian approach of probability theory as extended logic \cite{Jaynes:2003, Cox:1961} by imposing a symmetry and treating probability theory as syntactically complete. 

\section{Logic and Extended Logic}

To understand the completion I am proposing, we must understand certain aspects of logic that I propose to impose on probability theory.\\

Consider the following logical argument:
\begin{eqnarray*}
& A &\\
& B &\\[-9pt]
&\xdash[3em]& \\[-5pt]
\therefore & C &
\end{eqnarray*}
This is not a valid argument. The basic propositions ‘$A$’ and ‘$B$’ may have meanings for us that are not reflected in the formulation of the argument. For example, we may want the correspondence
\begin{align*}
A:&\hspace{2mm} \textrm{Socrates is a man.}\\
B:&\hspace{2mm} \textrm{If Socrates is a man, then Socrates is mortal.}\\
C:&\hspace{2mm} \textrm{Socrates is mortal.}
\end{align*}
A better formulation of the argument will then be to use logical implication $A \rightarrow C$ instead of $B$. We then get the new\footnote{Note here we are not considering uppercase propositions to be propositional variables; within each argument, propositions are constant. It is the arguments that change.} argument,
\begin{eqnarray*}
& A &\\
& A \rightarrow C & \\[-9pt]
&\xdash[3em] \\[-5pt]
\therefore & C &
\end{eqnarray*}
which is a valid argument. One may see from this a trivial aspect of logic; the validity of an argument is dependent on the structure of the argument. One needs to sufficiently define the structure in order to make an argument that is appropriate. \\

More importantly, the structure of an argument is the only thing the validity is dependent upon. Any influence to the validity of a logical argument beyond the form of the stated argument is extra-logical. One may have extra-logical influences in two ways; the choice of rules that define the logic used may be changed such that one uses a different logic; and one may have some meaning for a proposition in mind that is not defined within the argument. For this second influence, we shall take the position that one has insufficiently represented ones premises and thus the argument is not well formulated. For the first influence, this is a legitimate endeavour. We shall however stick to propositional logic and see how far we can go.\\

An argument is also independent of the labels one uses for the propositions; what is important is logical structure.\\

The position we are taking is sometimes called a syntactic logical interpretation of probability theory. Many philosophers take the position that there is meaning for some propositions relevant to scientific inference and prediction that cannot be defined through logical structure alone. These issues and others are discussed in the remarks.\\

We start by following Cox \cite{Cox:1961} who derived the product and sum rules from basic desiderata to extend logic. In the system, our primitive objects are degrees of plausibility
\begin{equation*}
A|B,
\end{equation*}
of $A$ given $B$, that are equal to real numbers. Probabilities are positive, continuous, monotonic functions, $f$, of plausibilities
\begin{equation*}
P(A|B)=f(A|B).
\end{equation*}
The function, $f$, is chosen such that certainty corresponds to the number 1. This choice is arbitrary but it leads to the particularly simple forms of the product and sum rules we have given.\\

We shall consider a degree of plausibility as directly analogous to a logical argument. This means two things:
\begin{enumerate}
\item The value of a plausibility (analogous to the validity of an argument) is dependent on only the explicit logical structure. From now on, any non-basic proposition will be written as $Z[A_i]$ or $Z[A_1, ... , A_n]$ as opposed to $Z$. 
When calculating a prior for $Z[A_i]$, the product and sum rules will constrain it to be functionally related to probabilities of basic propositions. We may then isolate the probabilities that cannot be calculated by using only the product and sum rules, and the rules of Boolean algebra.
\item Directly related to the above, the value of a degree of plausibility will not depend on the labels used on basic propositions. This gives us a powerful trick that is fundamental to derivations \cite{Jaynes:2003} of indifference and transformation group methods. For example,
\begin{equation*}
P(A_1|X[A_1,A_2 ] )=P(A_2|X[A_2,A_1 ] ).
\end{equation*}
If $X[A_1,A_2 ]$ is permutation symmetric; i.e., $X[A_1,A_2 ]=X[A_2,A_1 ]$, then
\begin{equation*} 
P(A_1|X[A_1,A_2 ] )=P(A_2|X[A_1,A_2 ] ),
\end{equation*}
which gives us a non-trivial constraint. Relabelling can also be used for individual probabilities separately within a functional relationship. For example,
\begin{eqnarray*}
P(A_1+A_2|A_3 )&=&P(A_1|A_3 ) + P(A_2|A_3 ) - P(A_1 A_2|A_3 ) \\
          &=& 2P(A_1|A_3 ) - P(A_1 A_2|A_3 ).
\end{eqnarray*}
\end{enumerate}

\section{Exclusivity, Exhaustivity and Indifference}

In nearly every probability calculation exclusivity and exhaustivity for some set of possibilities are (usually implicitly) assumed. Exclusivity means that if you assume that $A_i$ is true, then you may infer that any $A_j$, in some set $\left\{A_i\right\}^n_{i=1}$ where $i\neq j$, must be false. Exhaustivity means that you assume at least one of the set $\left\{A_i\right\}^n_{i=1}$ is true. These assumptions are often parametrised by the conditions
\begin{eqnarray*}
P(A_i A_j|I[A_i ] ) &=& P(A_i|I[A_i ] ) \delta_{ij}; \\
P(\sum^n_{i=1}A_i |I[A_i ]) &=& 1.
\end{eqnarray*}
If one wants to write a probability in terms of a mixture of others in the normal way, then it is necessary to make these assumptions;
\begin{eqnarray*}
P(B|X[A_i,B] )&=&P(B\sum^n_{i=1}A_i |X[A_i,B])\\
&=&\sum^n_{i=1}P(BA_i|X[A_i,B] ) \\ 
&=&\sum^n_{i=1}P(A_i|X[A_i,B] )P(B|A_i X[A_i,B] ).
\end{eqnarray*}
If one doesn't assume exclusivity, then one gets extra terms in the above equation.\\

One may consider exclusivity and exhaustivity as the sufficient logical definition of `possibility' and `property' \cite{Harrigan:2010, Pusey:2012}. A property may be considered a coarse grained possibility; for example, any classical observable in Hamiltonian mechanics segments the space of possible states, defining a property. In particular, the energy of a system is a property with various values of the energy related to various sets of possible states. The meaning that differentiates different types of properties and possibilities is defined by other logical relationships one assumes. So energy is a classification differentiated by causal structure (i.e., the equations one uses); we assume that if a system has a certain value for its energy at a certain time and is isolated, then the system will have the same energy at a later time; we have logical correlation. \\

One problem is that the conditions for exclusivity and exhaustivity are not derived from the explicit form of our assumptions $I[A_i]$. The explicit form of $I[A_i]$ for the simple set $\left\{A_1,A_2,A_3\right\}$ is as follows:
\begin{equation*}
M_3 [A_i ]=A_1\bar{A}_2\bar{A}_3+\bar{A}_1 A_2\bar{A}_3+ \bar{A}_1\bar{A}_2A_3+\bar{A}_1\bar{A}_2\bar{A}_3
\end{equation*}
for exclusivity and 
\begin{eqnarray*}
X_3 [A_i ] &=& A_1 A_2 A_3 + A_1 A_2\bar{A}_3 + A_1 \bar{A}_2 A_3 \\
&& + \bar{A}_1 A_2 A_3 + A_1\bar{A}_2\bar{A}_3+\bar{A}_1 A_2\bar{A}_3 \\
&& +\bar{A}_1 \bar{A}_2A_3
\end{eqnarray*}
for exhaustivity. This gives us 
\begin{equation*}
I_3 [A_i ] = M_3 [A_i ] X_3 [A_i ] = A_1\bar{A}_2\bar{A}_3 + \bar{A}_1 A_2 \bar{A}_3 + \bar{A}_1 \bar{A}_2 A_3.
\end{equation*}
These functions have been written in a non-minimal way where the function is a sum of terms with each term being a product of $A_i$'s and $\bar{A}_i$'s and containing $A_i$'s for all $1 \leq i \leq n$ for some $n$. In this form - which is often called disjunctive normal form - every term is exclusive by definition. Every propositional function can be written in this form. Each function can then be associated with a subset of the power set of $\left\{A_i\right\}^n_{i=1}$, $\left\{A_i\right\}^p$, where each term is associated with an element of $\left\{A_i\right\}^p$. The sum of terms associated with $\left\{A_i\right\}^p$ are exhaustive by definition. This will give us great flexibility in calculation. \\

Let us now calculate
\begin{eqnarray*}
P(A_1|I_3 [A_i ] ) &=& \frac{P(A_1\bar{A}_2\bar{A}_3|)}{P(A_1\bar{A}_2\bar{A}_3|) + P(\bar{A}_1A_2\bar{A}_3|) + P(\bar{A}_1\bar{A}_2A_3|)} \\
&=& \frac{1}{3} \times \frac{P(A_1\bar{A}_2\bar{A}_3|)}{P(A_1\bar{A}_2\bar{A}_3|)} \\
&=& \frac{1}{3},
\end{eqnarray*}
where we have used relabelling. Probabilities of the form $P(. |)$ refer to probabilities with minimal assumptions. This will be defined later. Indifference can be seen as fundamentally combinatoric in nature, where the probability is directly related to the number of ways one may assign a single non-negated proposition in a product of propositions. We can generalise to any $\left\{A_i\right\}^n_{i=1}$ in a simple manner.\\

One aspect of the above derivation is that the probabilities cancel out such that they do not need to be calculated. This is suggestive of why indifference could be derived in the past \cite{Jaynes:2003} without needing to go beyond the basic sum and product rules.\\

The above derivation also shows us that exclusivity and exhaustivity are not just necessary but also sufficient in deriving indifference.

\section{Determining Generic Probabilities}

We define a working set $\left\{A_i\right\}^n_{i=1}$ as the set of all propositions we are working with for a particular probability. This set may be made arbitrarily large:
\begin{eqnarray*}
P(A|B) &=& P(A|B)P(C + \bar{C}|AB)\\
&=& P(A[C + \bar{C}]|B)\\
&=& P(A|[C + \bar{C}]B).
\end{eqnarray*}
From this one can see we may add an arbitrary number of tautologies to the premises. One may consider this arbitrariness an important criterion for probability theory; we want our probabilities to be stable under arbitrary additions of tautologies to our assumptions. It is interesting that the product and sum rules give this to us for free.\\

An important thing to note is that we are allowing propositions like $C$ in the conclusions that have no representation in the premises.\\

Probabilities with minimal assumptions may be written as
\begin{equation*}
P(Z[A_1,...,A_n]|) = P(Z[A_1,...,A_n]|Q_n[A_1,...,A_n]),
\end{equation*}
where $Q_n[A_i] = \prod^n_{i=1}(A_i + \bar{A}_i)$. We may thus consider probabilities $P(.|)$ as ones either assuming nothing or only tautologies.\\

Let us now turn our attention to a generic probability
\begin{equation*}
P(Z[A_i ]|Y[A_i ]).
\end{equation*}
We may write
\begin{equation*}
P(Z[A_i ]|Y[A_i ] )=\frac{P(Z[A_i ]Y[A_i ]|)}{P(Y[A_i ]|)}.
\end{equation*}
Both of these factors may be written as sums of terms of the form $P(A_1 ... A_i\bar{A}_{i+1} ... \bar{A}_n|)$. These terms may be decomposed using the product rule. At this point the derivations of Jaynes \cite{Jaynes:2003} and Cox \cite{Cox:1961} give us no further support. Jaynes appeared \cite{Jaynes:2003} (P35) to consider probability theory complete as he expected one to always have background knowledge to determine the terms. Here we are explicitly assuming only tautologies with probabilities $P(.|)$ and hence cannot rely on such background knowledge. Moreover, as we are aiming at generality, we do not want to rely on such background knowledge.\\

To determine the terms, consider the following symmetry: The validity of an argument is invariant under swapping a basic proposition $A$ with its negation $\bar{A}$ in both the premises and the conclusions. Moreover, the swapping symmetry $A \leftrightarrow \bar{A}$ is a symmetry of logical structure. It may be seen as directly related to the double negation rule; imposing the symmetry on a trivial identity gives us the rule: 
\begin{equation*} 
\bar{A}=\bar{A}\rightarrow A=\bar{\bar{A}}.
\end{equation*}
Consider also that a possible state of affairs may be referred to by either $A$ or $\bar{A}$, with both choices giving equal consequences for argumentation. The two propositions are defined in contrast to one another (their truth values are opposed) and are not distinguished within the system in any other way. This lack of distinguishing factors is made more apparent when possibility is seen as an explicit extra assumption; the proposition $\bar{A}$ does not by definition mean that a proposition from a set of possibilities, other than $A$, must be true. Such meaning comes from an assumption $I_n[A_i]$.\\

I assert that our degrees of plausibility must satisfy the symmetry in order to not introduce an extra-logical bias into our framework. This may be considered part of the desideratum of consistency used in \cite{Jaynes:2003}.\\

Consider notation $x^j_{k}=P(A|A_1...A_j\Bar{B}_1...\Bar{B}_k)$. From our symmetry, one may impose\footnote{This condition is imposed on the plausibilities but may be stated in terms of probabilities.}
\begin{equation}
x^{j+1}_{k}=x^{j}_{k+1} \hspace{8mm} \forall j, k \geq 0.
\end{equation}
We now prove a lemma: $\forall j, k \geq 0,$ if $x^{j}_{k}=P(A|),$ then $x^{j}_{k+1}=P(A|)$.\\

Proof: Assume $x^{j_0}_{k_0} = P(A|) $ for some $j_0,k_0 \geq 0$. Then
\begin{eqnarray*}
P(A|) &=& P(A[A_q + \bar{A}_q]|A_1 ... A_{j_0} \bar{B}_1 ... \bar{B}_{k_0})\\
&=& P(A|)P(A|A_1 ... A_{j_0}A_q \bar{B}_1 ... \bar{B}_{k_0}) \\
&+& \left\{1 - P(A|)\right\}P(A|A_1 ... A_{j_0} \bar{B}_1 ... \bar{B}_{k_0} \bar{A}_q)\\
&=& P(A|)x^{j_0+1}_{k_0} + \left\{1 - P(A|)\right\}x^{j_0}_{k_0+1}\\
\end{eqnarray*}
Then by eq.(1), $x^{j_0}_{k_0+1}=P(A|).\hspace{5mm}\blacksquare$\\

From our lemma and (1), we get by induction from $x^0_0=P(A|)$,
\begin{equation*} 
\forall j, k \geq 0, x^{j}_{k}=P(A|).
\end{equation*}
Now see
\begin{eqnarray*}
P(A_1...A_j \bar{B}_1 ... \bar{B}_k|) &=& \prod^j_{l=1} \prod^k_{r=1}P(A_l|\mu_l)\left\{1 - P(B_r|\nu_r)\right\} \\
&=& P^j(A|)\left\{1-P(A|)\right\}^k,
\end{eqnarray*}
where $\mu_l$ and $\nu_r$ are products of various $A$'s and $B$'s. To determine $P(A|)$, we impose the symmetry again:
\begin{equation*}
A|=\bar{A}|.
\end{equation*}
From this condition one arrives at 
\begin{equation*}
P(A|)=\frac{1}{2}.
\end{equation*}
Our generic probability then becomes
\begin{equation*}
P(Z[A_i ]|Y[A_i ] )=\frac{M}{N}2^{n-m},
\end{equation*}
where $M$ and $N$ are just the number of terms in $Z[A_i ]Y[A_i ]$ and $Y[A_i ]$ respectively when the functions are written in minimal disjunctive normal form. Moreover, $M$ can heuristically be thought of as an unnormalised overlap between $Z[A_i ]$ and $Y[A_i ]$. The numbers $m$ and $n$ are the numbers of basic propositions in the terms of $Z[A_i ]Y[A_i ]$ and $Y[A_i ]$ respectively, written in minimal disjunctive normal form.

\section{Applications}

The precision and generality of our scientific statements are directly related to the precision and generality of the language used to make the statements. With the formulation of probability theory I have just proposed, we are able to determine precise probabilistic statements for a greater variety of situations then we were able to before. Immediate applications include situations where we do not or cannot assume exhaustivity and exclusivity:
\begin{enumerate}
\item We can calculate probabilities for propositions $A_i$ when we know only that they are one of m exhaustive possibilities;
\begin{equation*}
P(A_i|X_m [A_i ] )=\frac{2^{m-1}}{2^m-1}.
\end{equation*}
\item Consider a situation similar to one presented by Walley \cite{Walley:1996}: We have a bag of marbles. Suppose we know that they are labelled in a distinguishable way. In particular, they are numbered and we know there is a marble that is labelled with `1'. We know nothing about the number of marbles in the bag (perhaps the bag is magical, with the ability to hold an unlimited number of marbles). We want to know the probability that if we pick a marble from the bag, that marble will be the one labelled with `1'. This will generally depend on our knowledge of how we pick the marble. We are not interested in this particular aspect and if we know our method of picking cannot discriminate the labels, we may neglect this knowledge for our current purposes. 

Walley and others have proposed solutions to problems of this sort which go beyond the Bayesian framework. One sought after property of a probability in this situation is called regrouping invariance; i.e., it should somehow be invariant to changes in the `size of the sample space'. This presupposes that our probabilities are defined in terms of `sample spaces'.

Within the framework just proposed the solution requires only properly stating the salient assumptions; we have positive knowledge that there is a set of exclusive and exhaustive possibilities, we just do not know the size of the set. An appropriate probability will then be of the form
\begin{equation*}
\lim \limits_{n \to \infty}⁡P(A_1|\sum^n_{j=1}I_j [A_i]).
\end{equation*}
Note, assuming $\sum^n_{j=1}I_j [A_i]$ does not assume the various sample spaces are exclusive to each other. Exclusivity of sample spaces would require additional assumptions and change the probability. This is just one example of the precise choices we could make in our assumptions, exemplifying the generality of our approach.
\item Quantum theory has severe ontological problems. Our difficulty in solving these problems may be an insufficient formulation of probability theory \cite{Fuchs:2010}. Most if not all no-go theorems for ontological models of quantum theory \cite{Bell:1987, Harrigan:2010, Spekkens:2005, Hardy:2004, Pusey:2012} implicitly assume exclusivity and exhaustivity for the space of ontological states. The framework presented here allows for a whole class of models, which do not assume exclusivity and exhaustivity, to be explored.
\end{enumerate}

\section{Remarks}

The probabilistic framework here is considered as a symbolic system rather than a system of functions or measures on a predefined set. The framework is general enough to deal with situations where sets of possibilities are not assumed. The principle of indifference is derived as a consequence of our ability to relabel and the explication of the assumptions we implicitly make to define possibility. Indifference is thus not a principle imposed a priori or arbitrarily.\\

Probability theory as extended logic is completed by imposing a symmetry from propositional logic.\\

The degree to which one is convinced by the framework proposed here partly depends on whether one is convinced that propositional logic is sufficient for the task of scientific inference. We have seen how one may represent basic notions of possibility and property while still maintaining logical consistency. What propositional logic does not do are universals. I argue that universals are not directly relevant for scientific inference; a scientist would never be able to test the statement `all ravens are black'. \\

I propose the notion of universality is related to notions of induction and simplicity.\\

The framework just proposed does not directly justify induction. This is a good thing. An approach \cite{Carnap:1950} by Carnap - that has similar motivations to the approach here - tries to build induction directly into the framework. One problem is that the inductive predictions do not take into account ones assumptions; whether or not one predicts a sequence to continue at all and precisely how one predicts this depends on ones assumptions. Moreover, I submit these things should only depend on ones assumptions; if you make no assumptions you have no reason to predict the continuation of a sequence.\\

One may still perform inductive reasoning given certain assumptions such as a constant causal mechanism. There is, however, still a problem of induction: One may make valid predictions based on assumptions but those assumptions may not necessarily be justified. \\

The Bayesian framework has some built in notion of simplicity \cite{Jaynes:2003}(Ch.20). Consider two sets of propositional functions we'll call models, $\Omega_m$ and $\Omega_{m+1}$, where $\Omega_m$ is parametrised by $m$ parameters and $\Omega_{m+1}$ by $m+1$ parameters. Suppose the $m+1$'th parameter is $\theta_{m+1}$ and the subset $\Omega_{m+1}|_{\theta_{m+1}=0}$ has a one to one correspondence with $\Omega_m$ where each element in both sets is identified with one that produces the same likelihood for some data $D[A_i]$. We may take $\Omega_m$ and $\Omega_{m+1}$ as compound models, i.e., models where the parameters are unknown. If the elements are exclusive for both $\Omega_m$ and $\Omega_{m+1}$ and the point in the parameter space that gives a maximum likelihood (for data $D[A_i]$) is near $\theta_{m+1}=0$ and sharply peaked, then the likelihood for the compound model of $\Omega_m$ will generally be greater than the likelihood for $\Omega_{m+1}$; a set of models that predicts the observations as well as another but with less parameters will generally be better.\\

One limitation to this is that one has already chosen the sets of models to consider in a certain way. This has partly to do with ones preferences; do I judge a model with various parameters on the best choice of values of those parameters or do I judge a model on the total parameter space given to me? The choice also has to do with the choice of using a mathematical framework in the first place. In principle there are an infinite number of propositional functions that one may use as a model that have no discernable or consistent pattern. Can the restriction to propositional functions with consistent patterns be justified? This question becomes manifest in the proposed framework where we do not rely on calculating things with respect to a predefined set of alternative models; we may ask where those alternatives come from and why.\\

Note that the framework presented here manifests a primitive notion of simplicity for propositional functions themselves. The probability of some $Z[A_i]$, given no assumptions, is proportional to $2^{-m}$ where $m$ is the minimum number of propositions required to write $Z[A_i]$ in disjunctive normal form. The smaller the value of $m$, the `simpler' $Z[A_i]$ is.\\

I speculate that a justification for induction and simplicity comes from an assumption, $J[A_i]$, that restricts the set of propositional functions one may use. This restriction could be justified by epistemological considerations. Models with consistent patterns may then emerge due to combinatoric reasons.\\

The concept of possibility that is outlined in this article is suggestive of how scientific concepts may be defined generally. Possibility is a pattern of propositions within a model. Crucially, this pattern is not unique; different models with different sizes for possibility spaces will use different patterns (e.g., $I_2 [A_i]$ and $I_3 [A_i]$ are different). Moreover, the pattern may be nested such that the different possibilities are propositional functions rather than basic propositions. Within this framework, the concept of possibility cannot be defined as a form of classification, in contrast to some other attempts at the definition of a concept \cite{Valiant:1984}. I speculate that concepts like possibility and property may instead be associated with algorithms.\\

Universality may be defined as a concept.\\

This definition of concept suggests a motivation for its use. Consider an agent with data and assumption $J[A_i]$. There will likely be an infinite set of models to consider. Calculation for decisions may be computationally intractable. The agent may choose some scheme that best approximates the inferences one would ideally achieve. This scheme could involve algorithms for generating models. It may be the case that the best algorithms come from collections of nested concepts we may call general hypotheses. These general hypotheses may not give unique results but rather generate propositional functions dependent on input. Some of these general hypotheses may be well parametrised by mathematics.\\

Further work is required.


\begin{thebibliography}{0}%
\makeatletter
\providecommand \@ifxundefined [1]{%
 \@ifx{#1\undefined}
}%
\providecommand \@ifnum [1]{%
 \ifnum #1\expandafter \@firstoftwo
 \else \expandafter \@secondoftwo
 \fi
}%
\providecommand \@ifx [1]{%
 \ifx #1\expandafter \@firstoftwo
 \else \expandafter \@secondoftwo
 \fi
}%
\providecommand \natexlab [1]{#1}%
\providecommand \enquote  [1]{``#1''}%
\providecommand \bibnamefont  [1]{#1}%
\providecommand \bibfnamefont [1]{#1}%
\providecommand \citenamefont [1]{#1}%
\providecommand \href@noop [0]{\@secondoftwo}%
\providecommand \href [0]{\begingroup \@sanitize@url \@href}%
\providecommand \@href[1]{\@@startlink{#1}\@@href}%
\providecommand \@@href[1]{\endgroup#1\@@endlink}%
\providecommand \@sanitize@url [0]{\catcode `\\12\catcode `\$12\catcode
  `\&12\catcode `\#12\catcode `\^12\catcode `\_12\catcode `\%12\relax}%
\providecommand \@@startlink[1]{}%
\providecommand \@@endlink[0]{}%
\providecommand \url  [0]{\begingroup\@sanitize@url \@url }%
\providecommand \@url [1]{\endgroup\@href {#1}{\urlprefix }}%
\providecommand \urlprefix  [0]{URL }%
\providecommand \Eprint [0]{\href }%
\providecommand \doibase [0]{http://dx.doi.org/}%
\providecommand \selectlanguage [0]{\@gobble}%
\providecommand \bibinfo  [0]{\@secondoftwo}%
\providecommand \bibfield  [0]{\@secondoftwo}%
\providecommand \translation [1]{[#1]}%
\providecommand \BibitemOpen [0]{}%
\providecommand \bibitemStop [0]{}%
\providecommand \bibitemNoStop [0]{.\EOS\space}%
\providecommand \EOS [0]{\spacefactor3000\relax}%
\providecommand \BibitemShut  [1]{\csname bibitem#1\endcsname}%
\let\auto@bib@innerbib\@empty
\end{thebibliography}%


\begin{thebibliography}{14}



\bibitem{Jaynes:2003}
E.~T.~Jaynes,
``Probability Theory: The Logic of Science,''
Cambridge University Press, Cambridge (2003).

\bibitem{Finetti}
B.~de Finetti,
``Theory of Probability: A Critical Introductory Treatment,''
John Wiley \& Sons Ltd, New York (1974-75).

\bibitem{Shore:1980}
J.~E.~Shore and R.~W.~Johnson,
IEEE transactions on information theory, 26 (1980).

\bibitem{Cox:1961}
R.~T.~Cox,
``The Algebra of Probable Inference,''
John Hopkins Press, Baltimore (1961).

\bibitem{Walley:1996}
P.~Walley,
J.\ R.\ Statist.\ Soc.\ B, 58:3-57 (1996).

\bibitem{Fuchs:2010}
C.~Fuchs,
eprint arXiv:quant-ph/1003.5209v1 (2010).

\bibitem{Bell:1987}
J.~S.~Bell,
``Speakable and Unspeakable in Quantum Mechanics,''
Cambridge Univ. Press (1987).

\bibitem{Harrigan:2010}
N.~Harrigan, R.~W.~Spekkens,
Found.\ Phys.\ 40:125-157 (2010).

\bibitem{Spekkens:2005}
R.~W.~Spekkens,
Phys.\ Rev.\ A, 71:052108 (2005).

\bibitem{Hardy:2004}
L.~Hardy,
Stud.\ Hist.\ Phil.\ Mod.\ Phys., 35:267-276 (2004).

\bibitem{Pusey:2012}
M.~F.~Pusey, J.~Barrett, T.~Rudolph,
Nature Phys., 8:475 (2012).

\bibitem{Carnap:1950}
R.~Carnap,
``Logical Foundations of Probability,''
University of Chicago Press (1950).

\bibitem{Valiant:1984}
L.~G.~Valiant,
Communications of the ACM, 27:1134 (1984).

\bibitem{Cox:1946}
R.~T.~Cox,
Ann.\ J.\ Phys., 14:1-13 (1946).

\end{thebibliography}
\end{document}